\documentclass[11pt]{article}                        
\usepackage{amsfonts}
\usepackage{theorem}  \theorembodyfont{\upshape} 

\makeatletter
\def\draft{\textheight=10.5truein \textwidth=7.5truein \parindent=8pt
           \voffset=-1truein \topmargin=0Truein
           \ifcase \@ptsize \hoffset=-1.5truein \or \hoffset=-1.35truein
                        \or \hoffset=-1.15truein \fi}
\def\quality{\textheight=240mm \textwidth=160mm \topmargin=0Truein
             \ifcase \@ptsize \hoffset=-23mm
                     \or \hoffset=-20mm \or \hoffset=-15mm \fi}
\makeatother\quality

\def\bline(#1,#2)(#3,#4)(#5){\put(#1,#2){\line(#3,#4){#5}}}  

\newcommand\mlbscale{1pt} 
\newif\iffigs\figstrue 

\def\bfig(#1,#2)#3#4{\begin{figure} \begin{center}
    \framebox{\setlength{\unitlength}{\mlbscale}
       \iffigs \begin{picture}(#1,#2) #3 \end{picture}
       \else \begin{picture}(60,10)(0,0)
                   \put(0,0){\framebox(60,10){Figure}} \end{picture} \fi}
    \end{center} \caption{#4} \end{figure}}

\def\Bfig(#1,#2)#3#4{\begin{figure} \begin{center}
    \setlength{\unitlength}{\mlbscale}
       \iffigs \begin{picture}(#1,#2) #3 \end{picture}
       \else \begin{picture}(60,10)(0,0)
                   \put(0,0){\framebox(60,10){Figure}} \end{picture} \fi
    \end{center} \caption{#4} \end{figure}}

\def\bpic(#1,#2)#3{\setlength{\unitlength}{\mlbscale}
    \begin{picture}(#1,#2) #3 \end{picture}}

\def\n{\noindent} \def\IZ{{\mathbb{Z}}}   

\def\function#1{\left\{\!\!\!\begin{array}{ll} #1 \end{array} \right.}
 \def\a{\!\!\!&\!\!\!\!&}
\def\beq#1#2{\begin{equation} \label{#1} #2 \end{equation}}
\def\beaq#1#2{\label{#1} \begin{eqnarray} #2 \end{eqnarray}}

\def\thname{Theorem}     \def\lmname{Lemma}      \def\prname{Proposition}
\def\dfname{Definition}  \def\crname{Corollary}  \def\rmname{Remark}

\newtheorem{theorem}{\thname}[section]   
\newtheorem{lemma}{\lmname}[section]     
\newtheorem{corollary}[lemma]{\crname}   

\newtheorem{dftn}{\dfname}[section]
\newtheorem{rmrk}[lemma]{\rmname}

             \catcode`@=11 \@addtoreset{equation}{section} \catcode`@=12

\def\proof{\smallskip \noindent {\bf Proof. \ }}
\newcommand\emptysquare{\ \ $\Box$}
\def\qed{\hfill\emptysquare\linebreak\smallskip\par}

\def\*#1{#1^*}    \def\map{T}  
\def\IR{{\mathbb{R}}}  
\def\noprint#1{}  \def\FD{{\rm FD}}
   
 \def\v{{\rm v}}    \def\V{{\bar{V}}}
         \def\w{w}          \def\Reg{{\rm Reg}}
   \def\t#1{\tilde#1}

\def\intp#1{\left\lfloor#1\right\rfloor} \def\intpm#1{\lceil#1\rceil}
\def\toas#1{\stackrel{#1}{\longrightarrow}}
\def\blim#1#2{\if #1+ \limsup_{#2} \else {\if #1- \liminf_{#2} \else
              \lim_{#2}\left(\begin{array}{l}\sup\\
              \inf\end{array}\right) \fi} \fi} 

    \def\*#1{{#1^*}}
\begin{document}

\title{Hysteresis phenomenon in deterministic traffic flows}
\author{Michael Blank\thanks{
        Russian Academy of Sci., Inst. for
        Information Transm. Problems,
        and Observatoire de la Cote d'Azur, ~
        e-mail: blank@iitp.ru}
        \thanks{This research has been partially supported
                by Russian Foundation for Fundamental Research
                and French Ministry of Education grants.}
       }
\date{April 11, 2005} 

\maketitle

\n{\bf Abstract}. We study phase transitions of a system of
particles on the one-dimensional integer lattice moving with
constant acceleration, with a collision law respecting slower
particles. This simple deterministic ``particle-hopping'' traffic
flow model being a straightforward generalization to the well
known Nagel-Schreckenberg model covers also a more recent
slow-to-start model as a special case. The model has two distinct
ergodic (unmixed) phases with two critical values. When traffic
density is below the lowest critical value, the steady state of
the model corresponds to the ``free-flowing'' (or ``gaseous'')
phase. When the density exceeds the second critical value the
model produces large, persistent, well-defined traffic jams, which
correspond to the ``jammed'' (or ``liquid'') phase. Between the
two critical values each of these phases may take place, which can
be interpreted as an ``overcooled gas'' phase when a small
perturbation can change drastically gas into liquid. Mathematical
analysis is accomplished in part by the exact derivation of the
life-time of individual traffic jams for a given configuration of
particles.

\bigskip%
\n{\bf Keywords}: dynamical system, traffic flow, phase
transition, attractor.

\bigskip%
\n{\bf AMS Subject Classification}: Primary 37A60; Secondary
37B15, 37A50, 60K35.

\section{Introduction}\label{s:intro}

By a {\em traffic flow} we shall mean a collection of particles
moving along a straight line according to their velocities, and
the law describing how those velocities are changing (i.e., the
acceleration or deceleration of particles) is called the {\em
traffic flow model}.
For a long time theoretical analysis of traffic flow phenomena
has been dominated by hydrodynamic models in analogy to the
dynamics of viscous fluid (see, e.g. \cite{LeV,LW,Zh} and further
references therein). The main problem with this approach is that
it gives no information about the behavior of individual vehicles.
Moreover, the discrete nature of the traffic flow has some
features (like traffic jams) which do not have exact counterparts
in the hydrodynamic interpretation. Additionally, considered as a
kind of viscous fluid, the traffic flow turns out to have very
peculiar properties. First, it is very much `compressible': the
distance between two consecutive particles being originally quite
far from each other may shrink to zero under the dynamics due to a
large jam ahead of them. On the other hand, the opposite property
does not hold: the distance between two consecutive particles can
be enlarged at most by a constant depending on the parameters of
the model but not on time. Mention also that those hydrodynamic
models are difficult to treat in computer simulations of large
networks, while it is hard to compare parameters of the models
with empirical investigations. To overcome these difficulties
cellular automata models have been invented in the beginning of
90s (see \cite{CSS, GG} for reviews).

\Bfig(150,50)
      {
       \bline(0,40)(1,0)(150)   \put(154,38){$t$}
       \put(10,40){\circle*{5}} \put(10,45){\vector(1,0){20}}
       \put(45,40){\circle*{5}} \put(45,45){\vector(1,0){40}}
       \bline(0,0)(1,0)(150)    \put(154,-2){$t+1$}
       \put(30,0){\circle*{5}} \put(30,5){\vector(1,0){30}}
       \put(85,0){\circle*{5}} \put(85,5){\vector(1,0){50}}
       \bezier{30}(30,0)(30,22)(30,45)
       \bezier{30}(85,0)(85,22)(85,45)
       \put(35,40){\circle*{2}}
       \put(32,28){$i$} \put(8,28){$i_-$} \put(43,28){$i_+$}
       \put(10,50){$\intp{x_{i_-}}$}  \put(50,50){$\intp{x_{i_+}}$}
       \put(28,-12){$i'_-$} \put(83,-12){$i'_+$}
       \put(35,10){$x_{i_-}+a$}  \put(88,10){$x_{i_+}+a$}
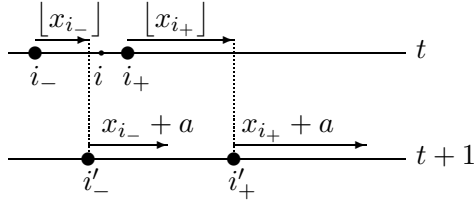
      }{Freely accelerating dynamics of particles. $i_\pm$ and
        $i'_\pm$ denote the positions of neighboring particles at
        moments of time $t$ and $t+1$ and $x_{\pm}$ -- corresponding
        velocities at time $t$.
        \label{f:free-motion}}

Consider a system of particles on the integer lattice $\IZ^1$
moving with constant acceleration $a\in\IR_+^1$ in a discrete time
with a collision law respecting slower particles. To make this
description precise we need a few definitions. Each particle is
described by the pair $(i,x_{i})$, where $i\in\IZ^{1}$ represents
the position of the particle and $x_{i}\in\IR^1$ its velocity,
which might be both positive and negative. Fix a configuration of
particles. For each $i\in\IZ^1$ by $i_-\le i$ we denote the site
containing the particle from our configuration closest to the site
$i$ from the left side, and by $i_+>i$ the one containing the
particle closest to $i$ from the right side (see
Fig.~\ref{f:free-motion}). Note the asymmetry of the definition:
$i_-$ might be equal to $i$, while $i_+$ is strictly larger than
$i$.

We shall say that a configuration of particles is {\em
admissible} if for any $i\in\IZ^{1}$ we have: %
\beq{e:adm}{i_{+}-i_{-}> \max\{0, ~ x_{i_{-}}, ~ -x_{i_{+}},
                                  ~ x_{i_{-}} - x_{i_{+}}\},
            \qquad |x_{i_\pm}|\le\v ,}%
where $\v>0$ is a parameter of the system describing the maximal
allowed velocity. In other words all particles in an admissible
configuration can be moved in arbitrary order by the distance
equal to the corresponding velocities without collisions.

On admissible configurations the dynamics is defined as follows.
First we change simultaneously the coordinates $(i,x_{i})$ of all
particles according to the rules: %
\beaq{e:qq}{i   \a\to  i + \intp{x_{i}} \label{e:dyn-adm1}\\
            x_i \a\to  \min\{x_{i} + a , \v \label{e:dyn-adm2}\} }%
where by $\intp{\cdot}$ we denote the integer part of a number,
i.e. $\intp{z}:=\max\{n\in\IZ: ~ n\le z\}$. Denote also by
$\intpm{\cdot}$ the smallest integer not smaller than the
considered number, i.e. $\intpm{z}:=\min\{n\in\IZ: ~ n\ge z\}$.

After this operation for each particle $(i,x_{i})$ for which the
admissibility condition (\ref{e:adm}) breaks down (which
correspond to a `collision' at the next step of the dynamics)
we correct its velocity: %
\beq{e:adm-cor1}{x_{i} \to \min\{x_{i}, ~ i_{+}-i_{-}-1\} .}%
Then for all particles where the velocities remain negative we
make an additional correction: %
\beq{e:adm-cor2}{x_{i} \to \max\{x_{i},
                 ~ i-(i-1)_{-}-1 - \max\{0,x_{(i-1)_{-}}\} \}. }%
It is straightforward to check that after these corrections the
configuration of particles becomes admissible.\footnote{We
   emphasize that the update rules (\ref{e:dyn-adm1}) --
   (\ref{e:adm-cor2}) are applied simultaneously to all
   particles in a given configuration.}

\begin{figure}
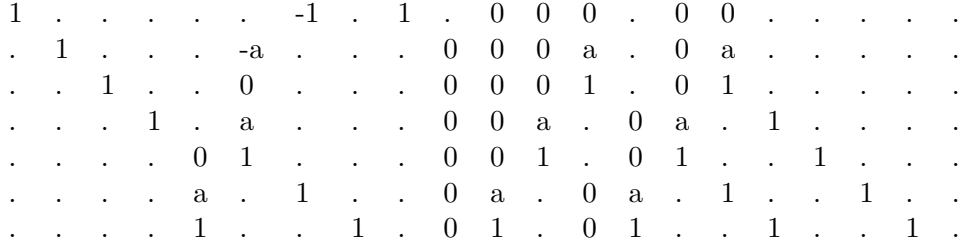
 \begin{center}
\begin{tabular}{llllllllllllllllllllll}
1&.&.&.&.&.&-1&.&1&.&0&0&0&.&0&0&.&.&.&.&.& \\
.&1&.&.&.&-a&.&.&.&0&0&0&a&.&0&a&.&.&.&.&.& \\
.&.&1&.&.&0&.&.&.&0&0&0&1&.&0&1&.&.&.&.&.& \\
.&.&.&1&.&a&.&.&.&0&0&a&.&0&a&.&1&.&.&.&.& \\
.&.&.&.&0&1&.&.&.&0&0&1&.&0&1&.&.&1&.&.&.& \\
.&.&.&.&a&.&1&.&.&0&a&.&0&a&.&1&.&.&1&.&.& \\
.&.&.&.&1&.&.&1&.&0&1&.&0&1&.&.&1&.&.&1&.& \\
\end{tabular}
\end{center}%
\vskip-0.5cm \caption{An example of the dynamics with $a=\frac12,
~ \v=1$. The positions of particles are marked by their
velocities and the positions of holes by dots.
\label{ex:model-dyn}}
\end{figure}%

Fig.~\ref{ex:model-dyn} illustrates the dynamics of our model
in the case $a=1/2$ and $\v=1$. Here and in the sequel we mark
the positions of particles in the configuration by their
velocities and the positions of holes by dots. The first line in
the figure corresponds to a finite segment of length 21 in the
initial configuration (i.e. at time 0) and the subsequent lines
describe the same segment for the moments of time from 1 to 6.
Note that we have chosen the segment such that the particles not
shown in the figure do not influence the dynamics of the ones in
the segment under consideration during the first six time steps.

Due to the constant acceleration after a finite time (of order
$\v/a$) all particles will start moving in the same (positive)
direction. Therefore, since we shall be interested mainly in the
limit behavior of the system as time goes to infinity, to simplify
notation we shall assume from now on that velocities of all
particles are nonnegative\footnote{Note, however, that due to the
   argument above all results in the paper remain valid in the
   initial setting with particle velocities of both signs.}  %
and shall describe the system by a {\em configuration}
$x=\{x_i\}_i\in\IR^\infty$, where a nonnegative entry $x_i$
corresponds to a particle at the site $i$ with the velocity
$x_i\ge0$, while all other sites (containing `holes') are set to
$-1$ (having in mind that the `holes' are moving in the opposite
direction to the particles). Under this assumption a configuration
$x$ is {\em admissible} if for any $i\in\IZ^{1}$ we have %
\beq{e:adm-pos}{i_{-} + x_{i_{-}}<i_{+}, \quad 
                x_i\in[0,\v]\cup\{-1\} .}%
The space of admissible configurations we denote by $X$. Using
this notation the dynamics can be defined as a superposition of
two maps $\map{}:={\cal A}\circ\sigma$, where the maps
$\sigma:X\to\IR^{\infty}$ and ${\cal A}:\IR^{\infty}\to X$ are
defined as follows: %
\beaq{e:maps}{%
 (\sigma x)_i \a:=
  \function{x_{i_-} &\mbox{if } i_-+\intp{x_{i_-}}=i\\
            -1      &\mbox{otherwise}} %
  , \\ %
 ({\cal A}x)_i \a:=
   \function{\min\{x_i+a,~i_+-i-1,~\v\} &\mbox{if } x_i\ge0\\
                   -1                   &\mbox{otherwise} } .} %
The map $\sigma$ is a (nonuniform) shift map, describing the
simultaneous free shift of all particles in the configuration $x$
by distances equal to their velocities, while the second map
${\cal A}$ describes the process of the acceleration/deceleration
of particles (see Fig.~\ref{f:free-motion}, ~\ref{ex:model-dyn}).
Observe that after the application of the map ${\cal A}$ any
configuration becomes admissible which ensures that the map $\map$
preserves the space of admissible configurations $X$.

The dynamics described by the map $\map$ for integer values of the
parameter $a>0$ coincides with the deterministic version of the
well known Nagel-Schreckenberg model (see e.g., \cite{CSS, GG}
for reviews). Note also that a version of the Nagel-Schreckenberg
model for continuous values of the velocities and spatial
coordinates was discussed in \cite{KWG1}.

We allow the velocity of a particle to be any real number between
$0$ and $\v$, but since the particles are moving on the integer
lattice the actual shift of a particle (described by the map
$\sigma$) is an integer, which means that a particle starts
moving only when its velocity becomes greater or equal to $1$.

\Bfig(100,100)
      {\footnotesize{
       \bline(0,0)(1,0)(160)   \bline(0,0)(0,1)(100)
       \bline(0,100)(1,0)(160) \bline(160,0)(0,1)(100)
       \put(20,20){\vector(1,0){120}} \put(20,20){\vector(0,1){70}}
       \thicklines
       \bline(20,75)(1,0)(26)  \bezier{200}(35,75)(70,23)(125,20)
       \thinlines
       \bezier{25}(46,75)(46,47)(46,20)
       \put(145,18){$\rho$}    \put(20,10){$\frac{1}{1+\w\v}$}
       \put(123,12){1}
       \put(12,85){$\V$}        \put(13,72){$\v$}
       \put(65,45){$\frac{\frac1\rho-1}{\w}$}
       \put(45,10){$\frac{1}{1+\v}$}
       \bezier{30}(35,75)(35,50)(35,20)
      }}
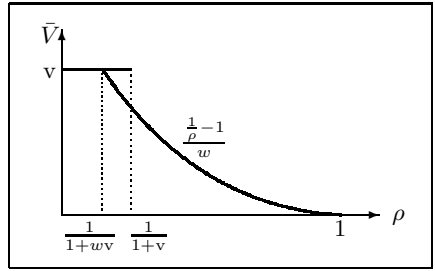
{Fundamental diagram for $\map{}$: dependence of the limit average
velocity $\V$ on the density of particles $\rho$. Here
$\w:=\intpm{1/a}$. \label{av-speed-den}}

One of the most striking observations in the theory of traffic
flows is the fact that in a large number of models as time goes
to infinity the limit `average velocity' of the flow depends only
on its density and does not depend on the finer characteristics.
This dependence is called `fundamental diagram' in physical
literature. Originally this observation has been established
numerically for finite lattices with periodic boundary
conditions, or in our terms for space periodic
configurations.\footnote{
  Indeed, an admissible configuration with a spatial period $\ell$
  remains periodic in space with the same spatial period $\ell$
  which explains the connection to numerical results.} %
Later this result has been proven rigorously for a broad class of
models \cite{Bl-jam, Bl-erg, NT}, assuming that the `average
velocity' is, indeed, the average velocity along a configuration,
which we shall call the `space average velocity'.

Recently a number of nonmarkovian traffic flow models
demonstrating the so called `metastable states' (see, e.g.
\cite{NFS, BHSS, KWG2} and references therein) have been
introduced. Numerous numerical simulations of those models show
that in distinction to the previous models the `space average
velocity' wildly fluctuates in time. Our own numerical results
also show the same phenomenon with fluctuations up to 100\% of the
`average' value if the acceleration $a<1$. To demonstrate this
choose a space periodic configuration with the spatial period of
length 6 consisting of the periodically repeating pattern
$[01..1.]$. Considering the dynamics with $a=\frac12, ~ \v=1$ of
only the main period of the configuration we get: %
$01..1. \to %
 a.1..0 \to %
 1..1.0 \to %
 .1..0a \to \dots$. %
As we see starting from the 2nd iterate the `space average
velocity' fluctuates in time periodically. In general the limit as
time goes to infinity `space average velocity' may not exist at
all, moreover, both upper and lower limits may fluctuate in time
as well. This shows that such approach is not adequate here. On
the other hand, the `space average velocity' tells almost nothing
about the movement of individual particles. Therefore it looks
reasonable to consider the averaging in time for individual
particles instead of the averaging in space. Note that due to the
nonergodicity of the process under study these two quantities need
not to coincide. On the other hand, it will be shown that the
`time average velocity' behaves much better and the main result of
this paper is the proof that the limit `time average velocity'
$\V(x)$ for an admissible regular\footnote{See definition in 
  Section~\ref{s:FD}.} %
configuration $x\in X$ having a density $\rho$ is described by the
following multivalued function %
\beq{e:00}{
    \FD_{a,\v}^{V}(\rho)
       = \min\{\v, ~ \frac{1/\rho-1}{\intpm{1/a}}\}
       ~ \cup ~ 1_{[(1+\intpm{1/a}\v)^{-1},
                    ~ (1+\v)^{-1}]}(\rho)\cdot\v  ,}%
corresponding to the fundamental diagram shown in
Fig.~\ref{av-speed-den}. The last term in the above relation
describes the upper (unstable) branch of the fundamental diagram.
Details related to lower and upper densities/velocities and
non-regular configurations will be discussed in
Section~\ref{s:FD}. Let me emphasize that non-regular 
configurations are not excluded and results related to the 
fundamental diagram are obtained for all configurations having 
densities.

Let us mention connections between the model under study and some
previously considered cases. In the case when all velocities of
particles are positive integers and the acceleration $a$ is an
integer as well we immediately recover the deterministic version
of the well known Nagel-Schreckenberg model introduced in
\cite{NS}. If one assumes that $a\ge\v$ then our system coincides
with the one studied analytically in \cite{Bl-jam, Bl-erg}. It is
known that in the case $a\ge\v$ the asymptotic behavior is much
simpler, in particular, only free flowing particles or free
flowing holes may show up, and, moreover, traffic jams cannot grow
in time. There is also another more technical difference due to
the fact that in that case the dynamics of holes is completely
symmetric to the dynamics of particles. This symmetry was heavily
used in \cite{Bl-jam, Bl-erg} for the analysis of the dynamics of
high density particle configurations. If $a<1$ there is no
independent dynamics of holes as such: to define their dynamics
one inevitably needs to take into account the velocities of 
particles. Therefore the analysis of a very complicated dynamics
of high density particle configurations cannot be reduced to the
low density case. Fortunately, as we shall see in the `jammed'
phase the dynamics of individual holes becomes periodic in time
(but not in space) after some initial stage. We shall see that if
$a<1$ then when the particle density exceeds some critical value
the traffic jams can grow in time and, moreover, an arbitrary
small jam (say, consisting of a few particles) may become
arbitrary large with time.

It is worth note that when the reciprocal to the acceleration $a$
is an integer and the initial configuration has only integer
velocities the model we are discussing coincides with the recent
`slow-to-start' model, which was extensively studied on the
numerical level (see e.g. \cite{SS,KWG2,NFS} and references
therein), but no analytical results have been obtained so far.
This model (being non-markovian in the sense that the dynamics
depends on the history) turned out to be the first one
demonstrating large fluctuations in `space average' statistics or
`metastable states'.

In this paper we restrict the analysis to the case of slow
particles with $\v=1$, because the analysis of the fast particles
with $\v>1$ and/or several lanes uses a rather different
mathematical apparatus -- substitution dynamics. Moreover, in the
general fast particles case (when $\v>1$) additionally to (static)
traffic jams that we consider in this paper there are new types of
(dynamic) jams when all particles in a jam might have positive
velocities and are moving as a whole with the velocity strictly
less than $\v$. These new types of jams lead to additional
branches in the fundamental diagram and will be discussed
elsewhere. Still since some of technical results obtained in the
paper remain valid for any $\v\in\IZ_{+}^{1}$ we shall keep the
notation $\v$ throughout the paper and shall specify explicitly if
only the case $\v=1$ is considered.

The paper is organized as follows. In the next Section we
introduce basic notions including the notion of the traffic jam 
and its characteristics such as a basin of attraction and a
life-time, and derive results describing their exact values as
functions of a given configuration. Using these results in
Section~\ref{s:FD} we prove the validity of the fundamental
diagram and provide its stability analysis in
Section~\ref{s:stability}. An alternative (more intuitive) proof
of the validity of the fundamental diagram via space-time
averaging is given in Appendix.

\section{Dynamics of traffic jams}\label{s:setup}

For a finite segment $x[n,m],
~~n<m\in\IZ_+^1\cup\{0\}=:\IZ_{\ge0}^1$ of a configuration $x\in
X$ we define the {\em density of particles} $\rho(x[n,m])$ as the
number of particles in the segment $N_{p}(x[n,m])$ divided by the
length of the segment, i.e. %
\beq{e:den-fin}{\rho(x[n,m]) %
    := \frac{N_{p}(x[n,m])}{m-n+1}
    = \frac1{m-n+1} \sum_{i=n}^m 1_{\ge0}(x_i) ,}%
where $1_{\ge0}(\cdot)$ is the indicator function of the set of
nonnegative numbers. The generalization of this notion for the
entire infinite configuration $x\in X$ leads to the notion of
{\em lower/upper densities}: %
\beq{e:den-ul}{\rho_{\pm}(x) :=
\blim{\pm}{n\to\infty}~\rho(x[1,n]), }%
where (and in the sequel) $\limsup$ corresponds to the index `+'
and $\liminf$ to the index `-'. If the lower and upper densities
coincide their common value $\rho(\cdot)$ will be called the {\em
density} of the configuration. According to Birkhoff ergodic
theorem for any translationally invariant measure on the integer
lattice the set of configurations having densities is the set of
full measure. In particular, for any space periodic configuration
the density is well defined.\footnote{Thus in numerical
   simulations only configurations with well defined densities
   can be observed.} %
Probably the analysis of configurations with well defined
densities would suffice, but since some configurations with
differing lower and upper densities lead to interesting behavior
we shall treat them as well.

The reason to fix $m=1$ in the definition above is that for any
given integer $m$ the limit points as $n\to\infty$ of the sequence
$\{\rho(x[m,n])\}_n$ are the same and do not depend on $m$. The
asymmetry with respect to the left and right `tails' of the
configuration reflects the fact that roles of those `tails' are
rather different in the long term dynamics. In particular,
consider two symmetric (to each other) configurations: $y^{(0)}$
having holes at all negative sites and particles at all positive
ones, and $y^{(1)}$ having particles at all negative sites and
holes at all positive ones. Any particle in the configuration
$y^{(0)}$ has zero velocity and will never move; on the other
hand, any given particle in $y^{(1)}$ will eventually with time
get the largest possible velocity $\v$.

Now we are ready to define time and space average velocities.
Since the case when there are no particles in the configurations
does not make any sense from the dynamical point of view, from
now on we shall assume that any admissible configuration contain
at least a finite number of particles. Let $L(x,i,t)$ be the
distance covered during $t$ time steps by the particle initially 
(at $t=0$) located at the site $i_-$ of the configuration $x$. 
Denote %
\beaq{}{\V^{\rm time}(x,i,t) \a\equiv \V(x,i,t) := \frac1t L(x,i,t) \\
     \V_\pm(x,i) \a:= \blim{\pm}{t\to\infty}\V(x,i,t) \\
     \V^{\rm space}(x[-n,m])
       \a:= \frac1{N_{p}(x[-n,m])}
            \sum_{i=-n}^m 1_{\ge0}(x_i) \cdot \intp{x_{i}} \\
     \V^{\rm space}_\pm(x)
       \a:= \blim{\pm}{n,m\to\infty}~\V^{\rm space}(x[-n,m]) .}
In other words $\V^{\rm space}$ describes the average velocity
along all particles in a given configuration, while $\V$
corresponds to the average in time velocity of a given particle.
More precisely these quantities describe average {\em effective}
velocities since we take into account only the actual movement of
particles and thus drop fractional parts of their velocities.

As we shall show the time average statistics (in distinction to
the space average one), does converge to the limit as time goes
to infinity and the result does not depend on the initial site
$i$, i.e. $\V(x,i)\equiv\V(x)$.

To this end, let us show that the upper and lower densities are
invariant with respect to dynamics.

\begin{lemma}\label{l:inv-den} $\rho_\pm(\map{}x)=\rho_\pm(x)$
for any configuration $x\in X$.
\end{lemma}
\proof For any $n,m\in\IZ_{\ge0}$ we have %
\beq{eq:mass-cons}{|N_{p}(x[-n,m]) - N_{p}((\map x)[-n,m])|
                   \le 2 ,}
since during one iteration of the map at most one particle can
enter the interval of sites from $-n$ to $m$ (from behind) and at
most one particle can leave this interval.\footnote{If only
   nonnegative velocities are considered 1 (instead of 2) would
   suffice in the inequality~(\ref{eq:mass-cons}), but in order
   to be valid in the original setting with velocities of both
   signs we put 2 here.} %

By the definition of the lower density there is a sequence of
integers $n_{j}\toas{j\to\infty}\infty$ such that
$$ \frac{N_{p}(x[1,n_{j}])}{n_{j}} \toas{j\to\infty} \rho_{-}(x) .$$
On the other hand, due to the inequality~(\ref{eq:mass-cons}) it
follows that $\rho_{-}(x,)$ is a limit point for partial sums for
the configuration $\map{} x$. Therefore we need to show only
that this is, indeed, the lower limit. Assume, on the contrary,
that there is another limit point, call it $\xi$, for the partial
sums for $\map{} x$ such that $\xi<\rho_{-}(x)$. Doing the
same operations with the partial sums for $\map{} x$ converging
to $\xi$ we can show that this value is also a limit point for the
partial sums for the configuration $x$, and, hence, $\xi$ cannot be
smaller than $\rho_{-}(x)$.

The proof for the upper density follows from the similar argument.
\qed

Introduce a {\em weight} function: %
\beq{e:weight}{\w(z)
    :=\function{\intpm{(1-z)/a} &\mbox{if } z\ge0 \\
                0               &\mbox{otherwise} }.}%
Referring to this function we shall drop the dependence on $z$ if
$z=0$ to simplify notation, i.e. $\w(0)\equiv \w$.

\begin{lemma}\label{l:vel-lim}
Let $x\in X$ and let the limit $\V(x,i)$ be well defined.
Then for any $j\in\IZ^1$ the limit $\V(x,j)$ is also well
defined and coincides with $\V(x,i)$.
\end{lemma}
\proof Observe that for any $x\in X$ and any $i\in\IZ^1$ the
particles at the sites $i_-$ and $i_+$ (call them left and right
ones) are consequent. After each iteration of the map $\map{}$
the distance between these particles changes by the difference
between their velocities, which might take values between $0$ and
$\v$. Since the left particle can be slowed down only by the
right one, we see that for any moment of time $t$ the distance
between the particles can be enlarged at most by $\w\v$. Indeed,
the longest time spent on the same site by the left particle
while the right one already started moving cannot exceed $\w$. On
the other hand, the maximal distance the right particle can cover
during this time is $\w\v$. Thus %
$$ 0 \le (i_+ + L(x,i_+,t)) - (i_- + L(x,i_-,t))
     \le i_+ - i_- + \w\v $$
or
$$ i_- - i_+ \le L(x,i_+,t) - L(x,i_-,t) \le \w\v .$$
Dividing by $t$ and using the definition of the time average
velocity we get %
\beq{eq:vel-inc}{
 |\V(x,i,t) - \V(x,i_+,t)|
 \le \max\{\frac1t \w\v, ~~ \frac1t(i_+ - i_-)\}
 \toas{t\to\infty}0 .}%
Thus $\V(x,i)=\V(x,i_+)$ for any $i$. Using the same argument for
$(i-1)_-$ and $i_+$ instead of $i$ one extends this result to
neighboring particles, and repeating this argument to all
particles in the configuration. \qed

This shows that if the limit $\V(x,i)$ exists for some site $i$
then the time average velocity is well defined, but we still have
to check the existence of at least one `good' site. On the other
hand, as we have shown the average velocity do not depend on the
left tail of a configuration. Therefore it is reasonable to
introduce for a configuration $x\in X$ the notion of its {\em
class of equivalence} $\hat{x}$ which includes all configurations
from $X$ which coincide with $x$ starting from some site, i.e.
$x,y\in\hat{x}$ if $\exists k\in\IZ^{1}: ~ x_{i}=y_{i} ~ \forall
i\ge k$.

Clearly in the absence of obstacles all particles are moving
freely. Therefore to understand the dynamics we need to study the
motion of `jams' as the only possible obstacles to the free
motion of particles. We say that a segment $x[n,m]$ with $n\le m$
corresponds to a {\em jam} if %
\beq{e:jam-def1}{0\le x_{i}<1 ~~ \forall i\in\{n,\dots,m\},}
\beq{e:jam-def2}{x_{n-1}<0, ~ |x_{m+1}|\ge1.} %
In other words the jam $x[n,m]$ is a locally maximal collection
of consecutive particles having zero velocities with the possible
exception of the {\em leading} one (located at the site $m$)
whose velocity is strictly less than $1$. For example, in the
configuration
$$~~\dots(a)\dots(0)\dots1\dots(00)1\dots(000a).1\dots~~$$
the jams are marked by parentheses.

The number of particles and their positions in a jam may change
with time: leading particles are becoming free and some new
particles are joining the jam coming from behind. However, only
one such change at a time might happen, and, in particular, a jam
cannot split into several new jams. Therefore we can analyze how
a given jam changes with time and the main quantity of interest
for us here is the minimal number of iterations after which the
jam will cease to exist. Denote by $J(t)$ the segment
corresponding to the given jam at the moment $t$ (in this
notation $J(0)$ is the original jam). Then by the {\em life-time}
of the jam $J$ we shall mean %
\beq{e:life-def}{\tau(J):=\sup\{t: ~ |J(t)|>0, ~ t>0\} ,}%
where $|A|$ is the length of the segment $A$.

`Attracting' the preceding particles, a jam plays a role similar
to an attractor in dynamical systems theory. Therefore it is
reasonable to study it in a similar way and to introduce the
notion of the {\em basin of attraction} (notation BA$(J)$) of the
jam $J:=x[n,m]$, by which we mean the segment $x[k,m]$ with $k\le
n\le m$ of the configuration $x$ containing all particles which
will eventually at positive time join the jam or can be stopped
by a particle from the jam $J$ during the time
$\tau(J)+1$.\footnote{The last assumption excludes the only
   possibility for a particle to create a new jam when it does not
   join an existing jam $J$ but is stopped at the time $\tau(J)+1$
   by the last particle in it which just get the velocity 1 but did
   not start moving yet. Consider an example: $.1.a \to ..01$. Here
   a particle preceding the jam $J$ containing the only particle
   with velocity $a$ is stopped at the next moment of time (and
   creates a new jam) when the jam $J$ already ceased to exist.} %
Examples of basins of attraction and their dynamics are shown in
Fig~\ref{ex:BA-dyn}.

Let us introduce the {\em weight} $W(x[n,m])$ of a segment $x[n,m]$
with $n\le m$ as %
\beaq{eq:weight-seg}{W(x[n,m])
   \a:= \w N_{p}(x[n,m-1]) + \w(x_m) \nonumber\\%
   \a\equiv \w\cdot(m-n)\cdot\rho(x[n,m-1]) + \w(x_m).}
The first term in this expression gives the contribution from all
particles in the segment except for the leading one, which is
described by the second term.

\begin{theorem}\label{t:life-time} Let $\v=1, ~ a\le1$ and let
$J(0):=x[n,m]$ with $n\le m$ be a jam. Then its basin of
attraction BA$(x[n,m])$ is the minimal segment
$x[k,m],~k\le{n}\le{m}$ for which relations %
\beq{eq:life-time-h}{N_{h}(x[k,m])+1 = W(x[k,m]) ,}%
\beq{eq:life-cond}{x_{k-i} < 1 - ia, \quad  i=1,2 }%
hold true,\footnote{
   When there is no segment satisfying (\ref{eq:life-time-h}) and
   (\ref{eq:life-cond}) we set $k=-\infty$, i.e. the BA coincides
   with the entire `left tail'.
   Note that the relations~(\ref{eq:life-time-h}), (\ref{eq:life-cond})
   give only an implicit information about the BA, which differs
   significantly from the result in \cite{Bl-jam, Bl-erg} where
   the much simpler case $a=\v$ has been considered.} %
where $N_{h}(x[k,m])$ is the number of holes in the segment
$x[k,m]$. Moreover, %
\beq{e:life-time}{\tau(J)=W({\rm BA}(J)) ,}%
i.e. the life-time of the jam is equal to the weight of its basin
of attraction, and under the dynamics the BA of a jam is
transformed to the BA of the remaining part of this jam.
\end{theorem}

Relations~(\ref{eq:life-time-h}), (\ref{eq:life-cond}) give a very
simple algorithm to find the left boundary $k$ of a BA: we move
to the left until the number of holes in the segment $x[k,m]$
will become equal to the weight of this segment minus one. After
this we check the condition (\ref{eq:life-cond}) and continue the
procedure if it does not hold or stop otherwise.

\proof First, let us rewrite the relation~(\ref{eq:life-time-h})
in a more suitable way. From the definition of the weight
function it follows that for $k\le m$ we have %
\beq{eq:num-part}{N_{p}(x[k,m]) = \intpm{W(x[k,m])/\w} .}%
On the other hand, we have the trivial identity
$$ m-k+1 = N_{p}(x[k,m]) + N_{h}(x[k,m]) ,$$
which together with (\ref{eq:num-part}) yields a new relation
equivalent to (\ref{eq:life-time-h}): %
\beq{eq:life-time}{m-k+2 = W(x[k,m]) + \intpm{W(x[k,m])/\w} .}%

Consider now in detail the dynamics of a BA of a jam. Two examples
are shown in Fig.~\ref{ex:BA-dyn}, where $a=1/2$ and we denote as
usual holes by dots and mark particles by their velocities. In
both examples the jams consist of particles having initially zero
velocities, and the right panel demonstrates the reason of the
inequalities~(\ref{eq:life-cond}), which exclude too slow
particles that cannot join the jam. Observe that if one changes
the velocity of the first particle in the right panel from $0$
to, say $a$, then the left boundary of the BA goes to the left
since the particle immediately preceding the left bracket will be 
stopped by the last particle in the former jam at the moment of 
time when it gets velocity 1.

\begin{figure}\begin{center}
\begin{tabular}{lllllllllllllllll}
$a=\frac12$\qquad&.&.&[.&1&.&.&0]&\qquad\qquad&.&0&[.&.&.&0&0]& \\
$\v=1$           & &.&.&[.&1&.&a]&            & &a&.&[.&.&0&a]& \\
                 & & &.&.&[.&0]& &            & &1&.&.&[.&0]& & \\
                 & & & &.&.&[a]& &            & & &1&.&.&[a]& & \\
                 & & & & &.&1& &              & & & &1&.&1& & \\
\end{tabular}
\end{center}%
\caption{Dynamics of basins of attraction whose boundaries
         are marked by square brackets.\label{ex:BA-dyn}}
\end{figure}%

Let us show that a basin of attraction of a jam either consists
of a single site, or its first site is occupied by a hole (i.e.
$x_{k}=-1$). To do this, consider a function
$$ \phi(\ell) := W(x[m-\ell,m]) + \intpm{W(x[m-\ell,m])/\w},
   \qquad \ell\in\IR_+^1 .$$
According to the definition of the weight of a
segment~(\ref{eq:weight-seg}) we can write the left boundary of
the BA$(x[n+1,m])=x[k,m]$ as $k=m-\t\ell$, where $\t\ell$ is the
smallest nonnegative solution to the equation %
\beq{eq:weight-jump}{\ell+2 = \phi(\ell) .}%
The function $\phi(\ell)$ is piecewise constant with jumps of
amplitude $(\w+1)$ at integer points corresponding to the lattice
sites $i$ where $x_i\ge0$ (i.e. where there is a particle in the
configuration $x$), and $\phi(0)=\w(x_m)+1$. Note that at a jump
point the value of $\phi$ corresponds to the upper end of the jump.

There might be two possibilities:
\begin{itemize}
\item $x_m\in[1-a,1)$. Then $\phi(0)=\w(x_m)+1=(\w-1)+1=\w$. Hence
      the smallest nonnegative solution to the
      equation~(\ref{eq:weight-jump}) occurs at the origin. This
      case corresponds to the BA consisting of a single site.

\item $x_m\in[0,1-a)$. In this case $\phi(0)=\w(x_m)+1>(\w-1)+1=\w$.
      Thus the starting point $\phi(0)$ of the graph $\phi(\ell)$
      is higher than the starting point of the straight line
      $\ell+\w$. Hence the first intersection between these graphs
      occurs at a horizontal piece of the graph of $\phi(\ell)$ 
      and thus it corresponds to a hole.
\end{itemize}

It remains to take into account the
inequalities~(\ref{eq:life-cond}). If at an intersection point
they are satisfied, then the result is proven, otherwise, due to
the presence of particles in the segment of length $2$ immediately
preceding to the interval $[m-\ell,m]$, the graph $\phi(\ell)$
makes some positive jumps of amplitude $\w+1$ each. Therefore
considering the graph after these jumps we are in the same
position as in the second case considered above and can apply the
same argument about the intersection point. Checking again the
inequalities~(\ref{eq:life-cond}) and continuing the procedure
(if necessary) we are finishing this part of the proof. 

Thus if $k<m$ then the first site $k$ of $x[k,m]$ is occupied by a
hole (i.e. $x_k<0$), while the last site $m$ is always occupied
by a particle with the velocity $x_m<1$ (since $x[n,m]$ is a jam).
Therefore after the application of the map $\map{}$ the velocity
of the particle at the site $m$ increases by $a$ and either this
particle becomes free and leaves the jam (if $x_m+a\ge1$, which
yields that the new leading particle of the jam is located at the
site $m':=m-1$), or the velocity remains below the threshold $1$
and $m':=m$. In both cases the weight of the segment $x[k,m']$
decreases by $1$. Hence, the solution to the
equation~(\ref{eq:weight-jump}) with $m=m'$ occurs now at the
point $\t\ell':=\t\ell+1$ which yields $k':=k+1$, i.e. the left
boundary of the BA is shifted by 1 position to the right, and the
resulting segment $x[k',m']$ is again the BA of the remaining
part $x[n,m']$ of the jam.

Note also that once we have found the solution to the
system~(\ref{eq:life-time-h}), (\ref{eq:life-cond}) then for all
subsequent time steps the inequalities~(\ref{eq:life-cond}) will
be satisfied automatically, since the particles immediately
preceding the BA cannot outstrip the left boundary of the BA
which is moving with the velocity 1.

Now, since after each iteration of the map $\map{}$ the weight of
the BA decreases by $1$, we deduce that the original weight of
the BA is equal to the life-time of the jam. \qed

Let us give also a short explanation to the ``physical meaning''
of the formula~(\ref{eq:life-time}). As we just have shown after
each iteration the left boundary of the BA increases by 1, while
the right boundary decreases by 1 once per $\w$ iterations, and
thus the ``space-time shape'' of the BA represents a skewed to the
right triangle (see Fig.~\ref{ex:BA-dyn}). When eventually the BA
size vanishes, the length of the original BA calculated using
these two observations gives the relation~(\ref{eq:life-time}).

\bigskip

This result shows that if all jams in the configuration $x\in X$
have finite BAs (or finite life-times, which is equivalent) then
eventually with time all particles will become free.

\section{Validity of the fundamental diagram}\label{s:FD}

Let us consider the dynamics of jams in more detail.

\begin{lemma}\label{l:ub-free} Let $\v=1, ~ a\le1$, $y\in X$ and let
$\rho_+(y)<\gamma_1:=(1+\w\v)^{-1}$. Then there is a configuration
$x\in\hat{y}$ in which only finite life-time jams may be present.
\end{lemma}

\proof Choose a representative $x\in\hat{y}$ satisfying the
assumption that $\rho(x[-n,0])\toas{n\to\infty}\rho_+(y)$ and
assume contrary to the statement of the Lemma that there exists a
jam $J=x[n,m]$ with $n<m<\infty$ having an infinite BA. (We do not
exclude the case $n=-\infty$ here.) Then by (\ref{eq:life-time})
using the same construction as in the proof of
Theorem~\ref{t:life-time} for any
$k<m$ we get: %
$$ m-k < W(x[k,m]) + \intpm{W(x[k,m])/\w} .$$
On the other hand, by (\ref{eq:weight-seg})
$$ W(x[k,m]) \le \rho(x[k,m])\cdot(m-k+1)\cdot \w .$$
Hence,
$$ \rho(x[k,m])\cdot(m-k+1)\cdot \w
 + \intpm{\rho(x[k,m])\cdot(m-k+1)} > m-k .$$
Passing to the limit as $k\to\infty$, we come to a contradiction:
$$ \rho_+(x)\ge(\w+1)^{-1}\equiv\gamma_1>\rho_+(x) .$$ \qed

As we shall see the critical value $\gamma_1:=(\w\v+1)^{-1}$ gives
the upper bound for the existence of a ``stable'' free-flowing
phase. One might think that under the conditions of
Lemma~\ref{l:ub-free} all configurations with densities below the
critical value $\gamma_1$ cannot have infinite life-time jams. To
show that this is not the case consider the configuration
$y^{(1)}$ having particles only at sites with negative
coordinates. Then $\rho(y^{(1)})=0$, but the segment
$y^{(1)}[-\infty,-2]$ of this configuration corresponds to a jam,
and since it has an infinite number of particles, the life-time of
this jam is infinite as well.

\begin{lemma}\label{l:hole-per} Let $\v=1, ~ a<1, ~ x\in X$ and
let in the configuration $x$ there exist a jam $J(t)$ with the
infinite BA (i.e., $|J(t)|>0 ~ \forall t\ge0$). Then each hole in
$x$, located originally to the right of $J(0)$, starting from the
moment of time when this hole meets with the particle that was the
leading one in the jam $J(0)$ will start moving periodically by
one position to the left exactly once per $\w$ time steps.
\end{lemma}

An immediate extension of this result is the following
observation. Let us call the largest moment of time at which a
given hole meets a particle the life-time of this hole. Then
under the assumption of Lemma~\ref{l:hole-per} the motion of each
hole having an infinite life-time is eventually periodic with the
period $\w$.

\smallskip

\proof Observe that since there is at least one infinite life-time
jam, then each particle belonging to its BA (i.e. located to the
left of $J$) will eventually join the jam $J$. On the other hand,
according to the dynamics the leading particle in the jam $J$ is
becoming free exactly once per $\w$ time steps. Thus these free
particles are passing holes (i.e. exchanging positions with them)
exactly once per $\w$ time steps (needed for the leading particle
in a jam to achieve the velocity 1 and to start moving).\footnote{
   Note that even when $\v>1$ and some holes may be trapped inside
   a jam, the situation remains the same due to a similar argument.
   Moreover, whence a particle will join some jam (non necessarily
   with an infinite BA) and will go through it, it will start
   passing holes exactly once per $\w$ time steps.}
\qed

This result provides us with the key observation to the calculation
of the time-average velocity in the jammed phase.

Assuming that a configuration $x\in X$ satisfies the assumptions
of Lemma~\ref{l:hole-per} let us consider in detail the dynamics
of a hole located initially at a site $j$ to the right of an
infinite life-time jam. According to Lemma~\ref{l:hole-per} after
some finite transient period this hole will start moving
periodically by one position to the left once per $\w$ time steps. 
Denote by $\vartheta_j$ the duration of the transient period and 
by $K_j$ the number of particles which the hole will meet during 
this time. Then we can characterize the deviation of the movement 
of the hole from the periodic one by the {\em defect} of the 
transient period $D_j:=K_j-\vartheta_j/\w$.

We shall say that a configuration $x\in X$ is {\em ultimately
jammed} if there is a representative $\hat x$ from the same
equivalence class in which for any $n\in\IZ^1$ there is a jam
$J_n(t)$ in the configuration $\hat x$ with the infinite BA and
$J_n(0)>n$, i.e. at the moment $t=0$ the jam $J_n(0)$ starts
further to the right from the site $n$. Since the dynamics of a
particle do not depend on particles located to the left from it we
shall assume (to simplify notation) that in an ultimately jammed
configuration all sites with nonpositive numbers are occupied by
particles with zero velocities.

Now we are ready to formulate the regularity property mentioned in
Introduction. Let in an ultimately jammed configuration $x$ holes
with nonnegative positions are located at sites $\{j_k\}$. We
shall say that the configuration $x$ is {\em regular} if the
functional %
\beq{e:def-reg-holes}{\Reg(x)
               := \limsup_{j_k\to\infty} |D_{j_k}|/j_k }%
describing statistics of absolute values of normalized defects
vanishes. In other words the regularity means that defects of
transient periods can grow only at a sublinear rate. Observe that
any spatially periodic ultimately jammed configuration is
regular\footnote{Clearly this property holds if distances between
   successive infinite life-time jams are uniformly bounded.
   See Lemma~\ref{l:reg-prop} for more general situations.} %
and that the property to be ultimately jammed is preserved under
the our equivalence relation.

The following statement describes the limit velocity statistics
for high density configurations.

\begin{lemma}\label{l:vel-jammmed}
Let $\v=1, ~ a<1, ~ x\in X$ and let $x$ be an ultimately jammed
configuration.\footnote{
    According to Lemma~\ref{l:ub-free} if the density is below
    $\gamma_1$ then all jams have finite BAs. Therefore the assumption
    about the presence of jams with infinite BAs yields the implicit
    assumption on the density.} %
Then %
\beq{eq:vel-av}{
  \left|\V_\pm(x) - \left(\frac1{\rho_\mp(x)} -1\right)\w^{-1}
  \right| \le \left(\frac1{\rho_\mp(x)} -1\right)\Reg(x) .}
\end{lemma}
\proof Observe first that we are in a position to apply
Lemma~\ref{l:hole-per}. For any $i\in\IZ^1$ the total distance
$L(x,i,t)$ covered by the particle initially located at the site
$i_-$ in the configuration $x\in X$ during the time $t>0$ is equal
to the number of holes encountered by the particle during this
time. For a given $t>0$ denote by $i_t$ the original position (at
moment $t=0$) in the configuration $x$ of the last hole our
particle meets during the time $t$. Then in the segment
$x[i_-,i_t]$ there are exactly $L(x,i,t)$ holes and, hence, its
length $i_t-i_-$ can be found
from the equation: %
$$ (i_t-i_-) \cdot(1 - \rho(x[i_-,i_t])) = L(x,i,t) ,$$
while the number of particles in this segment is equal to
$$ N_{p}(x[i_-,i_t]) = (i_t-i_-) \cdot\rho(x[i_-,i_t]) .$$
Thus %
\beq{e:part-def}{ N_{p}(x[i_-,i_t])
 = L(x,i,t) \cdot \frac{\rho(x[i_-,i_t])}{1 - \rho(x[i_-,i_t])} .}%

During the time $t$ the hole initially located at the site $i_t$
will meet with all those $N_{p}(x[i_-,i_t])$ particles. Since we
are interested in the long term dynamics we may assume that $t$
is so large that there is at least one infinite life time jam
between the sites $i_-$ and $i_t$, which guarantees that
$t>\vartheta_{i_t}$. Therefore %
$$ N_{p}(x[i_-,i_t]) = \frac{t-\vartheta_h(i_t)}{\w} + K_{i_t}
 = \frac{t}{\w} + D_{i_t} .$$
Substituting the obtained relation to equation~\ref{e:part-def} we
get
$$ \frac{t}{\w} + D_{i_t}
 = L(x,i,t) \cdot \frac{\rho(x[i_-,i_t])}{1 - \rho(x[i_-,i_t])} ,$$
which gives %
\beaq{e:vel-reg}{ \V(x,i,t) \a= \frac1t L(x,i,t) \nonumber\\
            \a= \frac1{\w}\left(\frac1{\rho(x[i_-,i_t])} - 1\right)
            + \left(\frac1{\rho(x[i_-,i_t])} - 1\right)
              \frac{D_{i_t}}t .}%
Observe now that $t\ge i_t$ and hence
$$ |D_{i_t}|/t \le |D_{i_t}|/i_t . $$
Therefore passing to the upper/lower limit as time goes to
infinity and taking into account that $i_t\toas{t\to\infty}\infty$
we get the result. \qed

Another more intuitive argument for the proof of this statement
based on a space-time averaging will be discussed in the Appendix.

\bigskip

The following result demonstrates the importance of the regularity
assumption for the convergence of the limit velocity.

\begin{lemma}\label{l:reg-necessity}
Let $x$ be an ultimately jammed configuration with $\Reg(x)>0$.
Then even if $\rho_-(x)=\rho_+(x)=\rho(x)<1$ the limit points of
the time average velocity may differ from the value
$(1/\rho(x)-1)/\w$.
\end{lemma}
\proof If the regularity assumption is not satisfied then there
exists a sequence of positions of holes $0<j_1<j_2<\dots$ such
that
$$ \lim_{k\to\infty} |D_{j_k}|/j_{k} = \Reg(x)>0 .$$
Assume contrary to our claim that $\V(x)=(1/\rho(x)-1)/\w$. Then
$i_t/t\toas{t\to\infty}\V(x)>0$ and therefore the second term in
the relation~(\ref{e:vel-reg}) describing fluctuations around the
time average velocity does not vanish with time. We came to the
contradiction. \qed

\begin{corollary} Let $x$ be an ultimately jammed configuration
with the density $\gamma_1<\rho(x)<\gamma_2$ and with $\Reg(x)>0$
and let the lower limit of normalized defects differ from the
upper one. Then using the same argument as above one can show that
the values $\V_\pm(x)$ differ as well.
\end{corollary}

\begin{lemma}\label{l:nec-in-many-jams}
The assumption in Lemma~\ref{l:vel-jammmed} that the configuration
$x$ is ultimately jammed is a necessary one and it is always
satisfied if %
\beq{ineq:vel-low}{\rho_-(x)>\gamma_2:=(1+\v)^{-1} .}%
\end{lemma}
\proof If there is only a finite number of jams then particles
originating from them might meet holes that were never met by
other particles before. Assume now that in the configuration $x$
only finite life-time jams can be present. Since a free particle
should have at least $\v$ holes immediately ahead of it and since
the length of the BA of a jam of $n$ particles exceeds $\w(n-1)$,
we deduce that in any segment of length $\ell$ there are at most
$\ell/(\v+1)$ particles, provided that $\ell$ is large enough.

It remains to show that under the assumption~(\ref{ineq:vel-low})
there is an infinite number of infinite life-time jams. Assume
that this is not the case, i.e. there exists $n\in\IZ^1$ such that
in the segment $x[n,\infty]$ only finite life-time jams are
present. Then we can apply the argument in the first part of the
proof to this segment to demonstrate that this assumption leads to
the inequality $\rho_-(x[n,\infty])<\gamma_2$. We came to the
contradiction since $\rho_-(x[n',\infty])=\rho_-(x[n,\infty])$ for
any $n'<n$. \qed

A more detailed analysis of the regularity functional improves
substantially the result of Lemma~\ref{l:vel-jammmed}. Moreover,
the sufficient condition for the regularity property can be
formulated in terms of {\em gaps} $G_k:=x[n_k,m_k]$ between
successive infinite life-time jams. Define %
\beq{e:def-reg-gaps}{\overline\Reg(x)
    :=\limsup_{k\to\infty}|G_k|/m_k .}%

\begin{lemma}\label{l:reg-prop} Let $\v=1, ~ a<1, ~ x\in X$ and
let $x$ be an ultimately jammed configuration. Then %
\begin{enumerate}
\item[(a)] $0\le\overline\Reg(x)\le1/2$, %
\item[(b)] $\overline\Reg(x)\ge\Reg(x)$,
           hence $\overline\Reg(x)=0$ implies $\Reg(x)=0$, %
\item[(c)] $\overline\Reg(x)=\Reg(x)\equiv0$ if
           $\rho_-(x)=\rho_+(x)>\gamma_2$.
\end{enumerate}
\end{lemma}
\proof Consider a gap $G_k$ between two successive infinite
life-time jams $J_k,J_{k+1}$. According to its definition the gap
may contain only free particles and finite life-time jams together
with their basins of attraction (otherwise those jams will have an
infinite life-time). Therefore for each particle in the gap there
should be at least one hole, which yields the assertion (a).

To prove the assertion (b) observe that for a hole located
initially at a site $j$ we have $0\le K_j\le\vartheta_j$ and thus
$$ |D_j| = |K_j-\vartheta_j/\w| \le (1-1/\w)\vartheta_j .$$
Now the statement follows from the observation that for a hole
located in the gap $G_k$ the length of the transient period cannot
exceed $|G_k|+\w$.

Assume now contrary to the last assertion that there is a
configuration $x\in X$ with $\rho_-(x)>\gamma_2$ and
$\overline\Reg(x)>0$. Then there is a sequence of gaps $G_k$ in
$x$ such that %
$$ \frac{|G_k|}{m_k}\toas{m_k\to\infty}\overline\Reg(x)>0 .$$ %
According to the proof of the assertion (a) the density of
particles in any gap $G_k$ cannot exceed $1/2$. Therefore
considering the configuration $x$ as a sequence of successive
alternative blocks of particles $J_k$ and gaps $G_k$ and taking
into account that lengths of jams $J_k$ should go to infinity with
$k$ (otherwise $\rho_-(x)\le\gamma_2=1/2$ since
$|G_k|\toas{k\to\infty}\infty$) we see that either
$\overline\Reg(x)=0$ or $\rho_-(x)<\rho_+(x)$ %
(since $x$ consists of alternative linearly growing blocks of two
types $J_i$ with $\rho(J_i)=1$ and $G_i$ with $\rho(G_i)\le1/2$).
\qed

From these technical statements we derive the main result about
the fundamental diagram.

\begin{theorem}\label{t:FD} Let $\v=1, ~ a\le1, ~ x\in X$. If %
\begin{itemize}
\item[(a)] $\rho_+(x)<\gamma_1:=(1+\w\v)^{-1}$ then $\V(x)=\v$,%
\item[(b)] $\rho_-(x)>\gamma_1$ and $x$ is not ultimately jammed
           (hence $\rho_-(x)<\gamma_2$) then $\V(x)=\v$,%
\item[(c)] $\rho_-(x)>\gamma_1$ and $x$ is ultimately jammed then
           $$\left|\V_\pm(x)
              - \left(\frac1{\rho_\mp(x)}-1\right)\w^{-1}\right|
            \le \left(\frac1{\rho_\mp(x)}-1\right)\Reg(x),$$%
\item[(d)] $\rho_-(x)=\rho_+(x)=\rho(x)>\gamma_2$ then %
           $\V(x)=(\frac1{\rho(x)}-1)\w^{-1}$.
\end{itemize}
\end{theorem}

\proof By Lemma~\ref{l:ub-free} under the condition
$\rho_{+}(x)<\gamma_{1}$ only jams with finite life-times may be
present in $x$. Consider a partition of $\IZ^{1}$ by
nonintersecting finite BAs of jams in the configuration $x$ and
their complement. Choose one of those BAs and denote by $i$ the
site containing the first particle preceding this BA. Even if
originally this particle does not have velocity one it will get
this velocity (and will start moving) at most after $\w$
iterations. On the other hand, by the definition of the BA this
particle (which does not belong to any BA to the right of it) will
never join a jam after the first moment of time when
$\V(x,i,t)=1$. Finally Lemma~\ref{l:vel-lim} finishes the proof of
item (a).

To prove the existence of the upper branch of the fundamental
diagram it is enough to show that for any
$0<\gamma<\gamma_2:=(1+\v)^{-1}$ there are configurations with
particle density $\gamma$ and such that eventually with time all
their particles will become free. To demonstrate this, consider a
space periodic configuration $x\in X$ with the space period of
length $2\ell$ in which only even sites
$0,2,\dots,\intp{2\ell\gamma}$ are occupied by particles with
velocity 1, while all other sites are occupied by holes. Clearly,
all particles in this configuration are free and will remain free
under dynamics. On the other hand, the density of this
configuration
$\rho(x)=\frac1{2\ell}\intp{2\ell\gamma}\toas{\ell\to\infty}\gamma$.

The remaining part of Theorem~\ref{t:FD} describing the lower
branch of the fundamental diagram follows immediately from
Lemmas~\ref{l:vel-jammmed},~\ref{l:nec-in-many-jams},~\ref{l:reg-prop}.
\qed

Obtained results give a complete characterization of the time
statistics of a configuration $x$ for which the density $\rho(x)$
is well defined or if it is not well defined but either
$\rho_+(x)<\gamma_1$ or $\rho_-(x)>\gamma_2$. Let us show that
there are other situations with much wilder behavior.

\begin{lemma}\label{l:01--10} Let $\v=1,~a\le1$. Then there exists
a configuration $y\in X$ such that
$\rho_-(y)=0, ~ \rho_+(y)=1$ and $\V_-(y)=0, ~ \V_+(y)=1$.
\end{lemma}
\proof We shall construct the configuration $y$ as follows. First
we set $y_i:=-1~~\forall i\le0$, i.e. we fill in the non positive
sites by holes. Then we shall fill the remaining sites by
alternative blocks $B_i$ consisting either only of particles with
zero velocities or only of holes. The lengths
$\ell_i:=|B_i|$ of those blocks we define inductively: %
$$ \ell_1:=1, ~~ \ell_2:=\w^{\ell_1}, ~~ \ell_3:=\w^{\ell_1+\ell_2},
   ~ \dots ~, ~  \ell_{k+1}:=\w^{\sum_{i=1}^k\ell_i}, ~\dots $$
Denoting $n_k:=\sum_{i=1}^k\ell_i$ we get $\ell_{k+1}=\w^{n_k}$
and $n_{k+1}=n_k + \ell_{k+1} = n_k + \w^{n_k}$. Therefore
$$ \rho(y[1,n_{2k+1}]) \ge \frac{\ell_{2k+1}}{n_{2k+1}}
 = \frac{\w^{n_k}}{n_k + \w^{n_k}} \toas{k\to\infty}1 .$$
Using a similar argument, but counting holes instead of particles
we get:
$$ 1 - \rho(y[1,n_{2k}]) \ge \frac{\ell_{2k}}{n_{2k}}
                             \toas{k\to\infty}1 .$$
Thus, $\rho(y[1,n_{2k}])\toas{k\to\infty}0$. This proves the
claim about the lower and upper densities.

Observe now that $\forall k\ge0$ the block $B_{k+1}$ corresponds to
a jam of length $\ell_{k+1}$ and that all particles to the left of
this jam lie in the segment $x[1,n_k]$, whose length is equal
to $\log_{\w}(\ell_{k+1})$. Therefore by Theorem~\ref{t:life-time}
all these particles belong to the basin of attraction of this jam
and hence will join it with time. On the other hand, since there
are no particles at negative sites, this shows also that life-times
of all jams in the configuration $y$ are finite. Note, however,
that the life-time of the $k$-th jam is of order $\exp(k)$.

Choose any particle of this configuration and consider how its
velocity changes in time. First from $t=0$ to $t=t_1\ge0$ the
particle might stay in a jam and have zero velocity. Then it is
becoming the leading one and during $\w$ time steps preserves its
position but accelerates until it get velocity $1$. After that the
particle starts moving freely with the velocity $1$ until it will 
catch up with the next jam. Then it will again stay in a jam 
having the zero velocity, etc. Due to the calculations above the 
duration of the alternative periods of staying in a jam and free 
moving (interrupted by short periods of acceleration of length 
$\w$) are growing exponentially. From this we get immediately the 
statement about the lower and upper time average velocities. \qed

\section{Stability/instability of the fundamental diagram}
\label{s:stability}

Qualitatively the main difference between the configurations
belonging to the upper and lower branches of the fundamental
diagram is related to the property to be ultimately jammed or not.
According to Theorem~\ref{t:FD} we need to study the stability of
this property only for configurations with densities in the region
$(\gamma_{1},\gamma_{2})$ where the two branches of the diagram
coexist. The following result demonstrate instability of the upper
branch of the fundamental diagram, while the stability of the
lower branch is discussed in Theorem~\ref{t:stability}.

\begin{theorem}\label{t:instability}
Let $\v=1, ~ a<1, ~ x\in X$ and let the density $\rho(x)$ be well
defined and $\rho(x)>\gamma_1$ while $\V(x)=1$ (i.e. the
configuration $x$ belongs to the upper branch of the fundamental
diagram). Then there exists an ultimately jammed configuration
$y\in X$ which differs from $x$ on a set $S$ having zero density.
\end{theorem}
\proof We shall construct the perturbed configuration $y$ as
follows. Set $\ell:=0$ and for for all $k\in\IZ_+^1$ for which
$x_k\ge0$ consider a sequence of numbers
$S_{k,\ell}:=\rho(x[\ell+1,k])$. We denote by $n_1$ the value of
$k$ which gives the first local maximum of $S_{k,\ell}$ for which
$S_{k,\ell}>\gamma_1$. Then we set $\ell:=n_1$ and continue the
procedure to find the number $n_2$ giving the local maximum
exceeding the value $\gamma_1$, etc. (see Fig.~\ref{f:averages}).
Eventually we shall have a monotonically growing sequence
$n_i\toas{i\to\infty}\infty$, from which we can further choose a
subsequence $\t{n}_i$ satisfying the assumption that $|\t{n}_i -
\t{n}_{i+1}|\ge2^i$. Observe that the set of integers
$\{\t{n}_i\}_i$ has zero density.

\Bfig(300,150)
      {\bline(0,0)(1,0)(300)   \bline(0,0)(0,1)(150)
       \bline(0,150)(1,0)(300) \bline(300,0)(0,1)(150)
       \bezier{100}(0,100)(150,100)(300,100)
       \bezier{50}(100,0)(100,75)(100,150)
       \bezier{50}(200,0)(200,75)(200,150)
       \thicklines
       \bezier{300}(0,0)(30,80)(50,30)
       \bezier{300}(50,30)(70,10)(100,140)
       \put(0,-8){$1$} \put(97,-8){$n_1$}
       \bezier{300}(100,0)(130,10)(150,60)
       \bezier{300}(150,60)(170,125)(200,130)
       \put(197,-8){$n_2$}
       \bezier{300}(200,0)(230,40)(250,30)
       \put(260,30){{\bf\dots}}
       \put(-10,1){$0$} \put(-10,97){$\gamma_1$}
       \put(-18,145){$S_{k,\ell}$} \put(290,-10){$\IZ^1$}
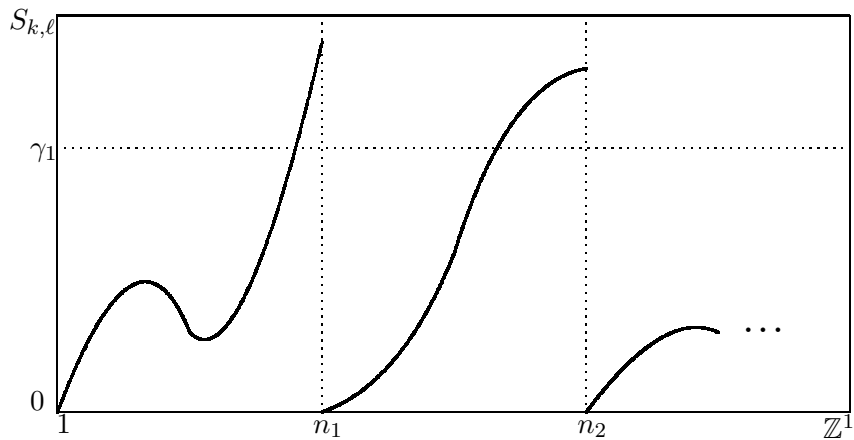
      }{Averages along a configuration \label{f:averages}}

Let $y_j:=0$ for all $j<=0$ and for positive $j$ we set
$y_j\equiv x_j$ except for the sites $\t{n}_i$ where we set
$y_{\t{n}_i}:=0$, i.e. we have changed the velocities of the
particles at the sites $\t{n}_i$ to 0 for all $i$. Since $\t{n}_i$
is the point of a local maximum, the site $\t{n}_i$ is the
leading point of some jam in $y$. Our claim is that the BA of
this jam covers the entire segment $[\t{n}_{i-1}+1,\t{n}_i]$.
Assume that this is not the case, i.e. there exists $N<\t{n}_i$
such that the BA corresponds to the segment
$y[\t{n}_i-N+1,\t{n}_i]$ and denote by $M$ the number of
particles in this segment. Due to the local monotonicity of the
sequence $S_{k,\ell}$ in the segment preceding the point
$\t{n}_i$ the strict inequality $\gamma_1<M/N$ holds true. On the
other hand, by Theorem~\ref{t:life-time} we have
$$ N = W(y[\t{n}_i-N+1,\t{n}_i])
       + \intp{W(y[\t{n}_i-N+1,\t{n}_i])/\w}
     = \w M + M .$$
Thus $M/N = \frac1{1+\w}=\gamma_1$, which contradicts to our
assumption.

Now, since the density of differing sites is equal to zero (due to
the choice of $\t{n}_i$) while the configuration $y$ satisfies the
assumptions of Lemma~\ref{l:vel-jammmed} we get the result. \qed

Clearly, changing back the zero density set of sites in the
ultimately jammed configuration $y$ constructed in the proof above
we get the configuration $x$ in which all particles will
eventually move freely. Therefore one might expect that all
ultimately jammed configurations with densities in the interval
$(\gamma_1,\gamma_2)$ are also unstable. The following statement
shows that this is not the case.

\begin{theorem}\label{t:stability} Let $\v=1, ~ a<1$. There exists
a regular configuration $x\in X$ with $\gamma_1<\rho(x)<\gamma_2$
and constants $0<A<B<\infty$ such that if a configuration $y\in X$
differs from $x$ at most at $A$ sites in any segment of length $B$
then $y$ is also regular. In other words, there is an open
neighborhood of the configuration $x$ consisting only of
configurations corresponding to the lower branch of the
fundamental diagram.
\end{theorem}
\proof Choose positive integers $N,M$ large enough such that
$\gamma_1<M/(N+M)<\gamma_2$ and consider a spatially periodic
configuration $x$ whose main period $x[1,N+M]$ consists of $N$
consecutive holes and $M$ consecutive particles having zero
velocities. Then the density $\rho=\frac{M}{N+M}$ of this
configuration is well defined and this configuration evidently
satisfies the conditions of Lemma~\ref{l:vel-jammmed}.

Our aim is to show that for a small enough (but still
nondegenerate) perturbation the property to have a countable
number of infinite life-time jams arbitrary far to the right is
preserved. Moreover we shall see that there exists a finite $N_1$
such that any segment of length $N_1$ in the perturbed
configuration $y$ intersects with some infinite life-time jam.
Clearly additional particles in $y$ (compared to $x$) should not
worry us and it is enough to consider only the case when under the
perturbation some particles are removed from the configuration $x$
or their velocities are changed.

Choose $n=k(N+M)$ with $k\in\IZ_+$ large enough and assume on the
contrary that all jams in $y[1,n]$ have finite life-times and
their BAs lie completely in this segment (i.e. they do not include
the site 0). Denote by $m_0$ the number of free particles in the
segment $y[1,n]$ and by $m_1$ the number of particles belonging to
jams in this segment.

By Theorem~\ref{t:life-time} we know that if $I$ is a finite
segment corresponding to a basin of attraction then
$W(I)>\w N_p(I)$ and $N_h(I)>W(I)$. Thus %
$$ |I| = N_p(I) + N_h(I) > (\w+1)N_p(I) .$$
Taking into account that each free particle is followed by a
hole we conclude that %
$$ 2m_0 + (\w+1)m_1 \le n .$$
On the other hand,
$$m_0+m_1=N_p(y[1,n])\le n\rho.$$
Thus %
$$ \frac{m_0}n \ge \frac{(\w+1)\rho -1}{\w-1} >0 .$$
Observe now that in order to get a free particle in the
configuration $y$ we need to change at least two sites in the
configuration $x$. From this one immediately can find the
constants $A,B$ setting %
$$ A:=\intp{\frac{(\w+1)\rho -1}{\w-1}\cdot B} $$
and choosing $B:=2k(N+M)$ to be so large that $A\ge1$. \qed

Using the same argument one can prove stability for a regular
configuration for which instead of space periodicity one assumes
that the gaps between infinite life-time jams are uniformly
bounded.

\smallskip

Let us show now that a single site perturbation of a configuration
consisting only of free particles can create a jam whose length
grows linearly with time. For $\v=1,~a=1/2$ consider a
configuration $x$ having free particles (with velocity 1) at all
even sites and holes at all others. We perturb this configuration
only at the origin setting the velocity at site 0 to zero and thus
creating a jam of length 1. According to Theorem~\ref{t:life-time}
the BA of this jam is infinite. On the other hand, after each
iteration a new particle joins this jam from the left while only
once per two iterations the leading particle leaves the jam, which
proves its linear growth.

\smallskip

A partial result in the direction of a measure-theoretic
interpretation of the notion of regularity for configurations in
the density region $(\gamma_1,\gamma_2)$ gives the following
Lemma.

\begin{lemma}\label{l:reg-density} Let a configuration $x$ has a
density $0<\rho(x)<1$. Decompose $x$ into alternative blocks of
particles $A_i$ and holes $B_i$: $x=\dots A_1B_1A_2B_2\dots$. Then
$$ \frac{|A_n|+|B_n|}{\sum_{i=1}^n(|A_i|+|B_i|)}
   \toas{n\to\infty}0 .$$ %
In other words, the lengths of blocks can grow only at a sublinear
speed.
\end{lemma}
\proof Let us prove that lengths of blocks of holes cannot grow at
a linear speed. Assume contrary to our claim that there exists a
constant $0<\zeta<1$ and a sequence of integers
$0<n_1<m_1<n_2<m_2<\dots$ with $n_i\le\zeta m_i$ for all $i$ such
that all sites in the segments $x[n_i,m_i]$ are filled by holes
while the segments $x[m_i,n_{i+1}]$ are filled by particles. Since
the density of the
configuration $x$ is well defined we have: %
$$ \rho(x[1,k])\toas{k\to\infty}\rho(x) .$$
On the other hand, %
$$ \rho(x[1,m_i]) = \frac{n_i\rho(x[1,n_i])}{m_i}
 \le \zeta \rho(x[1,n_i])\toas{i\to\infty}\zeta\rho(x) ,$$
which contradicts to the assumption that the density of $x$ is
well defined. A similar argument applies to blocks of particles as
well. \qed

Recalling the connection between the functionals $\Reg$ and
$\overline\Reg$ we see that the set of regular configurations do
not differ much from the set of all configurations having
densities even in the ``hysteresis'' region and we expect that for
any nontrivial translationally invariant measure on $X$ the set of
regular configurations is a set of full measure.

\section{Appendix: Proof of Lemma~\ref{l:vel-jammmed} via
         space-time averaging}\label{appendix}

Lemma~\ref{l:vel-jammmed} can be proven also in a more intuitive
way using a kind of a space-time averaging. Since this approach
explains the connection between the time and space averaging we
shall discuss it additionally to the formal proof given in
Section~\ref{s:FD}.

Observe that by the definition of the map $\map{}$ in the case
$\v=1$ we explicitly have the mass conservation law, i.e. after
each move one particle exchanges its position with one hole.
Therefore we can give a kind of `physical derivation' of the
relation~(\ref{eq:vel-av}). Assume that the (time) average
velocity $\V(x)$ is well defined and choose a segment of the
configuration $x$ of length $\ell$, provided $\ell\gg1$. Then the
number of particles in this segment is equal to $N_p:=\ell(\rho(x)
+ o(1/\ell))$, while the number of holes is equal to
$N_h:=\ell(1-\rho(x) + o(1/\ell))$. We are in a position to apply
Lemma~\ref{l:hole-per}, hence the holes are moving with the
average velocity $1/\w$ to the left. Thus in total $N_p$ 
particles will move to the distance $tN_p\V(x)$ during the time 
$t$, while in total $N_h$ holes will move to the distance 
$tN_h/\w$ during this time. Due to the mass conservation these 
quantities coincide and thus
$$ t\ell(\rho(x) + o(1/\ell))\V(x)
 = t\ell(1-\rho(x) + o(1/\ell))/\w .$$
Dividing by $t$ and $\ell$ and passing to the limit as
$\ell\to\infty$, we get
$$ \V(x) = \frac{1-\rho(x)}{\w \rho(x)} .$$

To make this argument rigorous, we need some additional work to be
done. We start with the space-periodic case and restrict ourselves
to just one period. For any given time interval the total shift of
all particles in the `spatial period' is well defined. Dividing
this total shift by the time and by the number of particles in the
`spatial period' and passing to the limit as time goes to infinity
(which exists due to the behavior of holes) we get the
`space-time' average velocity. On the other hand, since $\V(x,i)$
does not depend on $i$ we deduce that our `space-time' average
is, in fact, just $\V(x)$. Let us prove the last statement. Denote
the period length by $\ell$. For any moment of time $t$ we have
the following relation:
$$ \sum_{i=1}^{\ell}1_{\ge0}(x_{i})\cdot L(x,i,t)
 = \ell(1-\rho(x))\cdot \frac{t}{\w} + R(t) ,$$
Here (due to the usage of the indicator function) the summation is
taken over all particles in the `spatial period', the remaining
term $R(t)$ cannot exceed $\ell/\w$ on absolute value, and the
value $\ell(1-\rho(x))$ is equal to the number of holes in the
spatial period. Dividing by $t$, we get %
\beq{eq:av-vel-per}{
   \sum_{i=1}^{\ell}1_{\ge0}(x_{i})\cdot\V(x,i,t)
 = \ell(1-\rho(x))\cdot \frac{1}{\w} + O(1/t) ,}%
where $O(1/t)$ means a term of order $1/t$ as $t\to\infty$.

On the other hand, by the inequality~(\ref{eq:vel-inc}) for any
two particles initially located at sites $i,j$ of the `spatial
period' and any moment of time $t$ we have:
$$ |\V(x,i,t) - \V(x,j,t)|
  \le \ell\max\{\frac1t \w\v, ~~ \frac1t |i-j|\}
  \le \ell^2/t ,$$
provided that $\ell$ is large enough. Thus
$$ \left|\ell\rho(x)\cdot \V(x,\ell_-,t)
      - \sum_{i=1}^{\ell}1_{\ge0}(x_{i})\cdot \V(x,i,t)\right|
   \le \ell^3/t ,$$
where the value $\ell\rho(x)$ is equal to the number of particles
in the `spatial period'.

Therefore passing to the limit as $t\to\infty$ in the
relation~(\ref{eq:av-vel-per}) and using the inequality above, we
get
$$ \lim_{t\to\infty} \V(x,\ell_-,t)
 = \frac{\ell(1-\rho(x))}{\ell\rho(x)}\cdot \frac{1}{\w}
 = \left(\frac1{\rho(x)} -1\right)\w^{-1} ,$$
which finishes the proof in the space-periodic case.

\Bfig(250,60)
      {\bline(0,40)(1,0)(250)   \put(254,38){$t$}
       \bline(30,45)(1,0)(30)\bline(30,45)(0,-1)(5)
           \bline(60,45)(0,-1)(5) \put(35,48){$J_1(t)$}
       \bline(100,45)(1,0)(40)\bline(100,45)(0,-1)(5)
           \bline(140,45)(0,-1)(5)\put(105,48){$J_2(t)$}
       \bline(200,45)(1,0)(35)\bline(200,45)(0,-1)(5)
           \bline(235,45)(0,-1)(5) \put(205,48){$J_3(t)$}
       \put(60,35){\vector(1,0){80}}\put(70,35){\vector(-1,0){10}}
       \put(90,26){$I_1(t)$}
       \put(140,35){\vector(1,0){95}}\put(150,35){\vector(-1,0){10}}
       \put(170,26){$I_2(t)$}
       \bline(0,0)(1,0)(250)   \put(254,-2){$t'>t$}
       \bline(10,5)(1,0)(30)\bline(10,5)(0,-1)(5)
           \bline(40,5)(0,-1)(5) \put(15,9){$J_1(t')$}
       \bline(90,5)(1,0)(40)\bline(90,5)(0,-1)(5)
           \bline(130,5)(0,-1)(5)\put(95,9){$J_2(t')$}
       \bline(185,5)(1,0)(35)\bline(185,5)(0,-1)(5)
           \bline(220,5)(0,-1)(5) \put(190,8){$J_3(t')$}
       \put(40,-5){\vector(1,0){90}}\put(50,-5){\vector(-1,0){10}}
       \put(80,-15){$I_1(t')$}
       \put(130,-5){\vector(1,0){90}}\put(140,-5){\vector(-1,0){10}}
       \put(160,-15){$I_2(t')$}
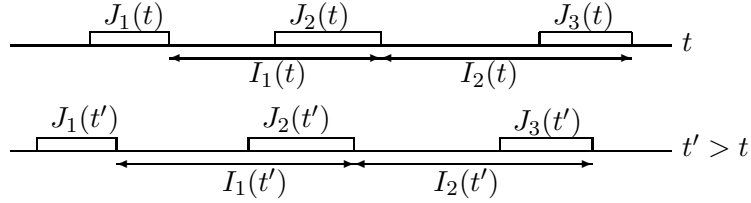
      }{Dynamics of traffic jams at moments of time $t$ and $t'>t$.
        \label{f:dyn-jams}}

In the general non space-periodic case one also can follow a
similar argument. According to our assumption in the
configuration $x$ there are finite jams $J_i(t)$ with infinite
life-times, denoted by rectangles in Fig.~\ref{f:dyn-jams}. Here
the parameter $t$ indicates the moment of time. The segment which
starts immediately after the site occupied by the leading
particle of the jam $J_i(t)$ and finishes at the site occupied by
the leading particle of the jam $J_{i+1}$, we denote by $I_i(t)$.
Then one can proceed as follows.

(a) The number of holes in $I_i(t)$ is invariant with respect to
the dynamics. Indeed, if a hole was located at time $t$ to the
right of a  jam $J(t)$, then for any moment of time $t'>t$ when
$J(t')$ still exists this hole should be still to the right of
$J(t')$. The reason is that the right boundary of a jam changes
only when its leading particle becomes free. However, the new
leading particle of the jam still remains to the left of our hole.

(b) The number of particles in $I_i(t)$ can change due to dynamics
at most by 1 in both sides. Indeed, the leading particle of each
the jams $J_i(t)$ is becoming free periodically with the period
$\w$, which is the only way how the number of particles in
$I_i(t)$ can change.

According to these properties one can think about the holes inside
of each segment $I_i(t)$ as a kind of a pump which is pushing the
particles through itself with the constant `productivity' equal to
$$ \frac{N_h(I_i(t))}{N_p(I_i(t))} \cdot \frac{1}{\w} ~,$$
which leads to the desired result.

\section*{Concluding remarks}

A deterministic generalization of the Nagel-Schreckenberg traffic
flow model (as well as of the slow-to-start model) with the real
valued acceleration has been proposed and studied analytically. It
has been shown that macroscopic physical properties of the model
under study such as the time average velocity of particles in the
flow depend crucially on the density of particles in the flow and
these results are described in terms of the corresponding
fundamental diagrams. Moreover, we have shown that instead of
space (or space-time) average velocities used earlier in traffic
models the time average velocity can be considered. It is proven
that this quantity characterizing the behavior of an individual
particle coincides for all particles in the flow and in
distinction to the space average velocity (which wildly
fluctuates in time) leads to a proper statistical description of
the dynamics.

One of the problems for future analysis is the question what
happens when the initial configuration is random or the dynamics
is weakly perturbed in random or deterministic sense (e.g., random
traffic lights are taken into account). The main idea in the paper
is that without perturbations one can predict the asymptotic
behavior of an individual configuration (and even of an individual
particle in it). Clearly integrating the results according to some
initial distribution one gets predictions for a family of random
initial configurations. One would expect that such questions might
be of interest only in the `hysteresis' region where both branches
of the fundamental diagram are present. Even there according to
Theorem 4.1 the influence of the upper branch should be
negligible. Another point is that it might be possible that under
arbitrary small perturbations the upper branch of the fundamental
diagram will disappear\footnote{From \cite{LZGM} it follows that
   if $\v=1,~a=1/2$ and random perturbations are applied only to
   jams but not to free particles the upper branch exists at
   least for spatially periodic configurations.} %
and that the non-regularity of initial configurations should not
matter in the random setting. However to check these predictions
one needs to study a probabilistic version of the model.

Finally let us mention that it is important to distinguish between
peculiarities related to a trivial divergence of various series
due to the fact that we deal with the infinite phase space
(infinite $\IZ^1$ lattice) and true finite-size phenomena taking
place even in the case of spatially periodic configurations (or a
finite lattice). The latter issue seems more important and our
result proves that the hysteresis phenomenon exists even in this
case.

The author is grateful to anonymous referees for very helpful
comments.

\end{document}